
%
%
%
%
%
%
\magnification=\magstephalf      
%
%
\vsize=7.5truein                 
\hsize=5.2truein                 
\newskip\stdskip                 
\stdskip=6pt plus3pt minus3pt    
\medskipamount=\stdskip          
\parindent=0pt                   
\parskip=\stdskip                
\abovedisplayskip=\stdskip       
\belowdisplayskip=\stdskip       
\mathsurround=0.75pt             
\overfullrule=0pt                
%
%
\def\ppar{\par\goodbreak\vskip 8pt plus 4pt minus 4pt}     
%
%
\def\stdspace{\hskip 0.75em plus 0.15em\ignorespaces}
\let\qua\stdspace 
%
%
%
%
%
%
%
\def\hexnumber#1{\ifcase#1 0\or 1\or 2\or 3\or 4\or 5\or 6\or 7\or 8\or
 9\or A\or B\or C\or D\or E\or F\fi}
%
%
\font\thirtnmsa=msam10 scaled 1315    
\font\tenmsa=msam10          \font\ninemsa=msam9
\font\sevenmsa=msam7         \font\sixmsa=msam6
\font\fivemsa=msam5
%
%
\newfam\msafam                  \textfont\msafam=\tenmsa
\scriptfont\msafam=\sevenmsa    \scriptscriptfont\msafam=\fivemsa
\edef\hexa{\hexnumber\msafam}        
\def\msa{\fam\msafam\tenmsa}         
%
%
\font\thirtnmsb=msbm10 scaled 1315   
\font\tenmsb=msbm10      \font\ninemsb=msbm9
\font\sevenmsb=msbm7     \font\sixmsb=msbm6
\font\fivemsb=msbm5
%
\newfam\msbfam                   \textfont\msbfam=\tenmsb       
\scriptfont\msbfam=\sevenmsb     \scriptscriptfont\msbfam=\fivemsb
\edef\hexb{\hexnumber\msbfam}    
\def\msb{\fam\msbfam\tenmsb}     
%
%
\font\thirtneufm=eufm10 scaled 1315   
\font\teneufm=eufm10                 \font\nineeufm=eufm9
\font\seveneufm=eufm7                \font\sixeufm=eufm6
\font\fiveeufm=eufm5
%
\newfam\eufmfam                    \textfont\eufmfam=\teneufm
\scriptfont\eufmfam=\seveneufm     \scriptscriptfont\eufmfam=\fiveeufm
\edef\hexf{\hexnumber\eufmfam}      
\def\frak{\fam\eufmfam\teneufm}     
%
%
%
\font\thirtnrm=cmr10 scaled 1315    
\font\ninerm=cmr9                   \font\sixrm=cmr6   
%
\font\thirtni=cmmi10 scaled 1315    
\font\ninei=cmmi9                   \font\sixi=cmmi6  
%
\font\thirtnsy=cmsy10 scaled 1315   
\font\ninesy=cmsy9                  \font\sixsy=cmsy6  
%
\font\thirtnbf=cmbx10 scaled 1315   
\font\ninebf=cmbx9                  \font\sixbf=cmbx6  
%
%
\font\thirtnex=cmex10 scaled 1315   
\font\nineex=cmex9                  
%
%
\font\thirtnit=cmti10 scaled 1315  
\font\nineit=cmti9                  
%
\font\thirtnsl=cmsl10 scaled 1315  
\font\ninesl=cmsl9                  
%
\font\thirtntt=cmtt10 scaled 1315  
\font\ninett=cmtt9                  
%
%
%
%
\def\small{%
%
%
\textfont0=\ninerm \scriptfont0=\sixrm \scriptscriptfont0=\fiverm
\def\rm{\fam0\ninerm}
%
%
\textfont1=\ninei \scriptfont1=\sixi \scriptscriptfont1=\fivei
%
%
\textfont2=\ninesy \scriptfont2=\sixsy \scriptscriptfont2=\fivesy
%
%
\textfont3=\nineex \scriptfont3=\nineex \scriptscriptfont3=\nineex
%
%
\textfont\bffam=\ninebf \scriptfont\bffam=\sixbf
\scriptscriptfont\bffam=\fivebf \def\bf{\fam\bffam\ninebf}%
%
%
\textfont\itfam=\nineit \def\it{\fam\itfam\nineit}%
\textfont\slfam=\ninesl \def\sl{\fam\slfam\ninesl}%
\textfont\ttfam=\ninett \def\tt{\fam\ttfam\ninett}%
%
%
%
\textfont\msafam=\ninemsa \scriptfont\msafam=\sixmsa
\scriptscriptfont\msafam=\fivemsa \def\msa{\fam\msafam\ninemsa}%
%
%
\textfont\msbfam=\ninemsb \scriptfont\msbfam=\sixmsb
\scriptscriptfont\msbfam=\fivemsb \def\msb{\fam\msbfam\ninemsb}%
%
%
\textfont\eufmfam=\nineeufm  \scriptfont\eufmfam=\sixeufm
\scriptscriptfont\eufmfam=\fiveeufm \def\frak{\fam\eufmfam\nineeufm}%
%
%
%
\normalbaselineskip=11pt%
\setbox\strutbox=\hbox{\vrule height8pt depth3pt width0pt}%
%
%
\normalbaselines\rm
%
%
\stdskip=4pt plus2pt minus2pt    
\medskipamount=\stdskip          
\parskip=\stdskip                
\abovedisplayskip=\stdskip       
\belowdisplayskip=\stdskip       
\def\ppar{\par\goodbreak\vskip 6pt plus 3pt minus 3pt}%
%
%
\def\section##1{\global\advance\sectionnumber by 1
\vskip-\lastskip\penalty-800\vskip 20pt plus10pt minus5pt 
\egroup{\bf\number\sectionnumber\quad##1}\bgroup\small         
\vskip 6pt plus3pt minus3pt
\nobreak\resultnumber=1}
}    
%
\def\beginsmall{\bgroup\small}
\let\endsmall\egroup
%
%
%
%
\def\large{%
\textfont0=\thirtnrm \scriptfont0=\ninerm \scriptscriptfont0=\sevenrm
\def\rm{\fam0\thirtnrm}%
\textfont1=\thirtni \scriptfont1=\ninei \scriptscriptfont1=\seveni
\textfont2=\thirtnsy \scriptfont2=\ninesy \scriptscriptfont2=\sevensy
\textfont3=\thirtnex \scriptfont3=\thirtnex \scriptscriptfont3=\thirtnex
\textfont\bffam=\thirtnbf \scriptfont\bffam=\ninebf
\scriptscriptfont\bffam=\sevenbf \def\bf{\fam\bffam\thirtnbf}%
\textfont\itfam=\thirtnit \def\it{\fam\itfam\thirtnit}%
\textfont\slfam=\thirtnsl \def\sl{\fam\slfam\thirtnsl}%
\textfont\ttfam=\thirtntt \def\tt{\fam\ttfam\thirtntt}%
\textfont\msafam=\thirtnmsa \scriptfont\msafam=\ninemsa
\scriptscriptfont\msafam=\sevenmsa \def\msa{\fam\msafam\thirtnmsa}%
\textfont\msbfam=\thirtnmsb \scriptfont\msbfam=\ninemsb
\scriptscriptfont\msbfam=\sevenmsb \def\msb{\fam\msbfam\thirtnmsb}%
\textfont\eufmfam=\thirtneufm  \scriptfont\eufmfam=\nineeufm
\scriptscriptfont\eufmfam=\seveneufm \def\frak{\fam\eufmfam\teneufm}%
\normalbaselineskip=16pt%
\setbox\strutbox=\hbox{\vrule height11.5pt depth4.5pt width0pt}%
\normalbaselines\rm}%
%
%
\def\Bbb#1{{\msb#1}}

%

\def\re{\Bbb R}
%
\mathchardef\plussquare="0\hexa01
\mathchardef\nge="3\hexb0B
\mathchardef\maltesecross="0\hexa7A
\mathchardef\del="0\hexf01
%
%
%
%
\font\sc=cmcsc10
%
%
%
%
\def\sqr#1#2{{\vcenter{\vbox{\hrule  height.#2truept
	\hbox{\vrule width.#2truept height#1truept 
	\kern#1truept \vrule width.#2truept}
	\hrule height.#2truept}}}}
\def\sq{\sqr55}    
%
%
%
%
\newcount\sectionnumber            
\newcount\resultnumber             
\sectionnumber=0\resultnumber=1    
%
%
%
\def\section#1{\global\advance\sectionnumber by 1
\xdef\nextkey{\number\sectionnumber}
\vskip-\lastskip\penalty-800\vskip 20pt plus10pt minus5pt 
{\large\bf\number\sectionnumber\quad#1}         
\vskip 8pt plus4pt minus4pt
\nobreak\resultnumber=1}      
%
%
%
%
%
         
%
%
%
%

%
\def\proc#1{\xdef\nextkey{\number\sectionnumber.\number\resultnumber}%
\vskip-\lastskip\ppar\bf%
\noindent#1\ \number\sectionnumber.\number\resultnumber
\stdspace\sl\global\advance\resultnumber by 1\ignorespaces}
 
%
%
\def\prf{\vskip-\lastskip\ppar\noindent{\bf Proof}%
\stdspace\rm}                            
\def\qed{\hfill$\sq$\par\goodbreak\rm}   
%
%
%
%
%
%
%
%
\def\proclaim#1{\vskip-\lastskip\ppar\bf%
\noindent#1\stdspace\sl\ignorespaces} 

%
%
%
%
\def\rk#1{\vskip-\lastskip\ppar{\bf #1}\stdspace\ignorespaces}                

%
%
%
%
%
%
\def\label{\xdef\nextkey{\number\sectionnumber.\number\resultnumber}%
\number\sectionnumber.\number\resultnumber
\global\advance\resultnumber by 1}
%
%
%
%
%
%
%
%
%
%
%
%
%
%
%
%
\newcount\refnumber              
\refnumber=1                     
\long\def\reflist#1\endreflist{%
\long\def\thereflist{#1}{\def\refkey##1##2\par{\xdef##1{\number\refnumber}%
\global\advance\refnumber by 1}%
\def\key##1##2\par{\expandafter\xdef%
\csname##1\endcsname{\number\refnumber}%
\global\advance\refnumber by 1}#1\par}}
\long\def\references{%
\penalty-800\vskip-\lastskip\vskip 15pt plus10pt minus5pt 
{\large\bf References}\ppar 
{\leftskip=25pt\frenchspacing    
\small\parskip=3pt plus2pt       
\def\refkey##1##2\par{\noindent  
\llap{[##1]\stdspace}\ignorespaces##2\par}         
\def\key##1##2\par{\noindent  
\llap{[\ref{##1}]\stdspace}\ignorespaces##2\par}  
\def\,{\thinspace}\thereflist\par}}
%
%
%
\newcount\footnotenumber         
\footnotenumber=1                
\def\fnote#1{\xdef\nextkey{\number\footnotenumber}%
{\small\ifnum\footnotenumber>9\parindent=14pt%
\else\parindent=10pt\fi\footnote{$^{\number\footnotenumber}$}%
{\hglue-5pt#1}\global\advance\footnotenumber by 1}}
%
%
%
%
%
%
%
\newcount\figurenumber          
\figurenumber=1                 
\def\caption#1{\xdef\nextkey{\number\figurenumber}%
\cl{\small Figure \number\figurenumber: #1}%
\global\advance\figurenumber by 1}
\def\figurelabel{\xdef\nextkey{\number\figurenumber}%
\cl{\small Figure \number\figurenumber}%
\global\advance\figurenumber by 1}
\long\def\figure#1\endfigure{{\xdef\nextkey{\number\figurenumber}%
\let\captiontext\relax\def\caption##1{\xdef\captiontext{##1}}%
\midinsert\cl{\ignorespaces#1\unskip\unskip\unskip\unskip}\vglue6pt\cl{\small 
Figure \number\figurenumber\ifx\captiontext\relax\else: \captiontext
\fi}\endinsert\global\advance\figurenumber by 1}}
%
%
%
%
%
%
%
\def\nextkey{??}   
%
\def\key#1{\expandafter\xdef\csname #1\endcsname{\nextkey}}
\def\ref#1{\expandafter\ifx\csname #1\endcsname\relax
\immediate\write16{Reference {#1} undefined}??\else
\csname #1\endcsname\fi}
%
%
%
%
%
%
%
\newread\gtinfile
\newwrite\gtreffile
\def\useforwardrefs{
\openin\gtinfile\jobname.ref
\ifeof\gtinfile
\closein\gtinfile
\immediate\write16{No file \jobname.ref}
\else
\closein\gtinfile
\input \jobname.ref
\fi
\immediate\openout\gtreffile \jobname.ref
%
%
\def\key##1{{\def\\{\noexpand}%
\expandafter\xdef\csname ##1\endcsname{\nextkey}%
\immediate\write\gtreffile{\\\expandafter\\\def\\\csname ##1\\\endcsname%
{\nextkey}}}}
%
%
\long\def\reflist##1\endreflist{%
\long\def\thereflist{##1}{\def\refkey####1####2\par{\xdef####1{%
\number\refnumber}{\def\\{\noexpand}\immediate\write\gtreffile
{\\\def\\####1{\number\refnumber}}}\global\advance\refnumber by 1}%
\def\key####1####2\par{\expandafter\xdef%
\csname####1\endcsname{\number\refnumber}%
{\def\\{\noexpand}\immediate\write\gtreffile
{\\\expandafter\\\def\\\csname ####1\\\endcsname{\number\refnumber}}}
\global\advance\refnumber by 1}##1\par}}
\long\def\biblio##1\endbiblio{\reflist##1\endreflist\references}%
%
%
\def\numkey##1{{\def\\{\noexpand}%
\xdef##1{\number\sectionnumber.\number\resultnumber}
\immediate\write\gtreffile{\\\def\\##1%
{\number\sectionnumber.\number\resultnumber}}}}
\def\seckey##1{{\def\\{\noexpand}\xdef##1{\number\sectionnumber}
\immediate\write\gtreffile{\\\def\\##1{\number\sectionnumber}}}}
\def\figkey##1{\xdef##1{\number\figurenumber}%
{\def\\{\noexpand}\immediate\write\gtreffile%
{\\\def\\##1{\number\figurenumber}}}
\number\figurenumber\global\advance\figurenumber by 1}
}   
%
%
%
%
\def\figkey#1{\xdef#1{\number\figurenumber}%
\number\figurenumber\global\advance\figurenumber by 1}
\def\fig#1#2\endfig{%
\midinsert\cl{#2}\vglue6pt\cl{\small Figure #1}\endinsert}
\def\newfig{\number\figurenumber\global\advance\figurenumber by 1}
\def\numkey#1{\xdef#1{\number\sectionnumber.\number\resultnumber}}
\def\seckey#1{\xdef#1{\number\sectionnumber}}
%
%
%
%
%
%
%
%
%
\def\verb{\catcode`\"=\active}       
\def\brev{\catcode`\"=12}            
\brev                                
\verb                                
{\obeyspaces\gdef {\ }}              
{\catcode`\`=\active\gdef`{\relax\lq}}
\def"{%
\begingroup\baselineskip=12pt\def\par{\leavevmode\endgraf}%
\tt\obeylines\obeyspaces\parskip=0pt\parindent=0pt%
\catcode`\$=12\catcode`\&=12\catcode`\^=12\catcode`\#=12%
\catcode`\_=12\catcode`\~=12%
\catcode`\{=12\catcode`\}=12\catcode`\%=12\catcode`\\=12%
\catcode`\`=\active\let"\endgroup}
\brev      
%
%
%
%
%
%
\def\item#1{\par\leavevmode\llap{#1\stdspace}%
\ignorespaces}                             
%
%

%
%
\def\co{\colon\thinspace}    
\def\np{\vfil\eject}         
\def\nl{\hfil\break}         
\def\cl{\centerline}         
\def\gt{{\mathsurround=0pt\it $\cal G\mskip-2mu$eometry \&\ 
$\cal T\!\!$opology}}        
\def\gtm{{\mathsurround=0pt\it $\cal G\mskip-2mu$eometry \&\ 
$\cal T\!\!$opology $\cal M\mskip-1mu$onographs}}    
\def\agt{{\mathsurround=0pt\it$\cal A\mskip-.7mu$lgebraic \&\ 
$\cal G\mskip-2mu$eometric $\cal T\!\!$opology}}  
%
%
%

%
%
%
%
%
\def\title#1{\def\thetitle{#1}}

\def\author#1{\edef\previousauthors{\theauthors}
 \ifx\theauthors\relax\def\theauthors{#1}\else
 \def\theauthors{\previousauthors\par#1}\fi}

%
\def\address#1{\edef\previousaddresses{\theaddress}
 \ifx\theaddress\relax\def\theaddress{#1}\else
 \def\theaddress{\previousaddresses\par\vskip 2pt\par#1}\fi}
\def\secondaddress#1{\edef\previousaddresses{\theaddress}
 \ifx\theaddress\relax\def\theaddress{#1}\else
 \def\theaddress{\previousaddresses\par{\rm and}\par#1}\fi}   

\def\email#1{\edef\previousemails{\theemail}
 \ifx\theemail\relax\def\theemail{#1}\else
 \def\theemail{\previousemails\hskip 0.75em\relax#1}\fi}
\def\secondemail#1{\edef\previousemails{\theemail}
 \ifx\theemail\relax\def\theemail{#1}\else
 \def\theemail{\previousemails\hskip 0.75em{\rm and}\hskip 0.75em
 \relax#1}\fi}
\def\url#1{\edef\previousurls{\theurl}
 \ifx\theurl\relax\def\theurl{#1}\else
 \def\theurl{\previousurls\hskip 0.75em\relax#1}\fi}
\def\secondurl#1{\edef\previousurls{\theurl}
 \ifx\theurl\relax\def\theurl{#1}\else
 \def\theurl{\previousurls\hskip 0.75em{\rm and}\hskip 0.75em
 \relax#1}\fi}
\long\def\abstract#1\endabstract{\long\def\theabstract{#1}}
\def\primaryclass#1{\def\theprimaryclass{#1}}
\def\secondaryclass#1{\def\thesecondaryclass{#1}}
\def\keywords#1{\def\thekeywords{#1}}
%
%
\let\\\par\let\thetitle\relax\let\theshorttitle\relax
\let\theauthors\relax\let\theshortauthors\relax
\let\theaddress\relax\let\theshortaddress\relax
\let\theemail\relax\let\theurl\relax
\let\theabstract\relax\let\theprimaryclass\relax
\let\thesecondaryclass\relax\let\thekeywords\relax
%
%
%
%
\long\def\maketitlepage{    

\vglue 0.2truein   

%
{\parskip=0pt\leftskip 0pt plus 1fil\def\\{\par\smallskip}{\large
\bf\thetitle}\par\medskip}   

\vglue 0.15truein 

%
{\parskip=0pt\leftskip 0pt plus 1fil\def\\{\par}{\sc\theauthors}
\par\medskip}%
 
\vglue 0.1truein 

%
{\small\parskip=0pt
{\leftskip 0pt plus 1fil\def\\{\par}{\sl\theaddress}\par}
\ifx\theemail\relax\else  
\vglue 5pt \def\\{\stdspace{\rm and}\stdspace} 
\cl{Email:\stdspace\tt\theemail}\fi
\ifx\theurl\relax\else    
\vglue 5pt \def\\{\stdspace{\rm and}\stdspace} 
\cl{URL:\stdspace\tt\theurl}\fi\par}

\vglue 7pt 

{\bf Abstract}

\vglue 5pt

\theabstract

\vglue 7pt 

{\bf AMS Classification numbers}\quad Primary:\quad \theprimaryclass\par

Secondary:\quad \thesecondaryclass

\vglue 5pt 

{\bf Keywords:}\quad \thekeywords

\np  

}    
%
%
\long\def\makeshorttitle{    


%
{\parskip=0pt\leftskip 0pt plus 1fil\def\\{\par\smallskip}{\large
\bf\thetitle}\par\medskip}   

\vglue 0.05truein 

%
{\parskip=0pt\leftskip 0pt plus 1fil\def\\{\par}{\sc\theauthors}
\par\medskip}%
 
\vglue 0.03truein 

%
{\small\parskip=0pt
{\leftskip 0pt plus 1fil\def\\{\par}{\sl\ifx\theshortaddress\relax
\theaddress\else\theshortaddress\fi}\par}
\ifx\theemail\relax\else  
\vglue 5pt \def\\{\stdspace{\rm and}\stdspace} 
\cl{Email:\stdspace\tt\theemail}\fi
\ifx\theurl\relax\else    
\vglue 5pt \def\\{\stdspace{\rm and}\stdspace} 
\cl{URL:\stdspace\tt\theurl}\fi\par}

\vglue 10pt 


{\small\leftskip 25pt\rightskip 25pt{\bf Abstract}\stdspace\theabstract

{\bf AMS Classification}\stdspace\theprimaryclass
\ifx\thesecondaryclass\relax\else; \thesecondaryclass\fi\par
{\bf Keywords}\stdspace \thekeywords\par}
\vglue 7pt
}    
\let\maketitle\makeshorttitle        
%
%

\def\volumenumber#1{\def\thevolumenumber{#1}}
\def\volumename#1{\def\thevolumename{#1}}
\def\volumeyear#1{\def\thevolumeyear{#1}}
\def\pagenumbers#1#2{\def\startpage{#1}\def\finishpage{#2}}
\def\published#1{\def\publishdate{#1}}

\volumenumber{X}
\volumename{Volume name goes here}
\volumeyear{20XX}
\pagenumbers{1}{XXX}
\published{XX Xxxember 20XX}

\long\def\makegtmontitle{   

\count0=\startpage

\gtm\nl        
{\small Volume \thevolumenumber: \thevolumename\nl 
Pages \startpage--\finishpage\nl}

\vglue 0.1truein   

{\parskip=0pt\leftskip 0pt plus 1fil\def\\{\par\smallskip}{\large
\bf\thetitle}\par\medskip}   
\vglue 0.05truein 

%
{\parskip=0pt\leftskip 0pt plus 1fil\def\\{\par}{\sc\theauthors}
\par\medskip}%
 
\vglue 0.03truein 


{\small\leftskip 25pt\rightskip 25pt{\bf Abstract}\stdspace\theabstract

{\bf AMS Classification}\stdspace\theprimaryclass
\ifx\thesecondaryclass\relax\else; \thesecondaryclass\fi\par
{\bf Keywords}\stdspace \thekeywords\par}\vglue 7pt

}   

\long\def\makeagttitle{   
\agt\hfill      
\hbox to 60truept{\vbox to 0pt{\vglue -14truept{\bf [Logo here]}\vss}\hss}
\break
{\small Volume \thevolumenumber\ (\thevolumeyear)
\startpage--\finishpage\nl
Published: \publishdate}

\vglue .2truein

{\parskip=0pt\leftskip 0pt plus 1fil\def\\{\par\smallskip}{\large
\bf\thetitle}\par\medskip}   
\vglue 0.05truein 

%
{\parskip=0pt\leftskip 0pt plus 1fil\def\\{\par}{\sc\theauthors}
\par\medskip}%
 
\vglue 0.03truein 


{\small\leftskip 25truept\rightskip 25truept{\bf Abstract}\stdspace\theabstract

{\bf AMS Classification}\stdspace\theprimaryclass
\ifx\thesecondaryclass\relax\else; \thesecondaryclass\fi\par
{\bf Keywords}\stdspace \thekeywords\par}\vglue 7truept

}   


\def\Addresses{\bigskip
{\small \parskip 0pt \leftskip 0pt \rightskip 0pt plus 1fil \def\\{\par}
\sl\theaddress\par\medskip \rm Email:\stdspace\tt\theemail\par
\ifx\theurl\relax\else\smallskip \rm URL:\stdspace\tt\theurl\par\fi}}

\def\agtart{
\hoffset 14truemm
\voffset 31truemm
\font\phead=cmsl9 scaled 950
\font\pnum=cmbx10 scaled 913
\font\pfoot=cmsl9 scaled 950
\headline{\vbox to 0pt{\vskip -4.5mm\line{\small\phead\ifnum
\count0=\startpage ISSN numbers are printed here
\hfill {\pnum\folio}\else\ifodd\count0\def\\{ }%
\ifx\theshorttitle\relax\thetitle\else\theshorttitle\fi\hfill{\pnum\folio}
\else\def\\{ and }{\pnum\folio}\hfill\ifx\theshortauthors\relax\theauthors
\else\theshortauthors\fi\fi\fi}\vss}}
\footline{\vbox to 0pt{\vglue 0mm\line{\small\pfoot\ifnum\count0=\startpage
Copyright declaration is printed here\hfill\else
\agt, Volume \thevolumenumber\ (\thevolumeyear)\hfill\fi}\vss}}
\let\maketitle\makeagttitle\let\makeshorttitle\makeagttitle
\let\maketitlepage\makeagttitle}

\def\gtmonart{
\hoffset 14truemm
\voffset 31truemm
\font\phead=cmsl9 scaled 950
\font\pnum=cmbx10 scaled 913
\font\pfoot=cmsl9 scaled 950
\headline{\vbox to 0pt{\vskip -4.5mm\line{\small\phead\ifnum
\count0=\startpage ISSN numbers are printed here
\hfill {\pnum\folio}\else\ifodd\count0\def\\{ }%
\ifx\theshorttitle\relax\thetitle\else\theshorttitle\fi\hfill{\pnum\folio}
\else\def\\{ and }{\pnum\folio}\hfill\ifx\theshortauthors\relax\theauthors
\else\theshortauthors\fi\fi\fi}\vss}}
\footline{\vbox to 0pt{\vglue 0mm\line{\small\pfoot\ifnum\count0=\startpage
Copyright declaration is printed here\hfill\else
\gtm, Volume \thevolumenumber\ (\thevolumeyear)\hfill\fi}\vss}}
\let\maketitle\makegtmontitle\let\makeshorttitle\makegtmontitle
\let\maketitlepage\makegtmontitle}

\def\gtart{
\hoffset 14truemm
\voffset 31truemm
\font\phead=cmsl9 scaled 950
\font\pnum=cmbx10 scaled 913
\font\pfoot=cmsl9 scaled 950
\headline{\vbox to 0pt{\vskip -4.5mm\line{\small\phead\ifnum
\count0=\startpage ISSN numbers are printed here
\hfill {\pnum\folio}\else\ifodd\count0\def\\{ }%
\ifx\theshorttitle\relax\thetitle\else\theshorttitle\fi\hfill{\pnum\folio}
\else\def\\{ and }{\pnum\folio}\hfill\ifx\theshortauthors\relax\theauthors
\else\theshortauthors\fi\fi\fi}\vss}}
\footline{\vbox to 0pt{\vglue 0mm\line{\small\pfoot\ifnum\count0=\startpage
Copyright declaration is printed here\hfill\else
\gt, Volume \thevolumenumber\ (\thevolumeyear)\hfill\fi}\vss}}
\let\maketitle\maketitlepage\let\makeshorttitle\maketitlepage}

\input epsf
%
%
\def\relabelbox{%
  \hbox\bgroup%
}%
\def\endrelabelbox{%
\egroup%
}%
\def\relabel #1#2 {%
  \special{ps:/a {} def}%
  \smash{\rlap{#2}}%
}%
\def\adjustrelabel <#1,#2> #3#4 {%
  \special{ps:/a {} def}%
  \smash{\rlap{\kern #1 \raise #2\hbox{#4}}}%
}%
\def\extralabel <#1,#2> #3 {\smash{\rlap{\kern #1 \raise #2\hbox{#3}}}}%

\chardef\newinsCatAt\the\catcode `\@
\catcode `\@=11
%
%
%
\newskip\insertskipamount\newskip\inserthardskipamount
\insertskipamount 12pt plus2pt  
\inserthardskipamount 4pt       
\def\insertskip{\vskip\insertskipamount}
%
%
\newskip\LastSkip
\def\SaveLastSkip{\LastSkip\lastskip}
\def\RestoreLastSkip{\nobreak\vskip-\LastSkip\vskip\LastSkip}
%
%
\newcount\SplitTest
\def\SetSplitTest{\SplitTest\insertpenalties
  \insert\topins{\floatingpenalty1}%
  \advance\SplitTest-\insertpenalties}
%
%
\def\midinsert{\par
 \SaveLastSkip\penalty-150\SetSplitTest\RestoreLastSkip
 \ifnum\SplitTest=-1
  \@midfalse\p@gefalse\else\@midtrue\fi\@ins}
\def\@ins{\par\begingroup\setbox\z@\vbox\bgroup%
  \vglue\inserthardskipamount}
\def\endinsert{\egroup 
  \if@mid \dimen@\ht\z@ \advance\dimen@\dp\z@
    \advance\dimen@\insertskipamount
    \advance\dimen@\pagetotal\advance\dimen@-\pageshrink
    \ifdim\dimen@>\pagegoal\@midfalse\p@gefalse\fi\fi
  \if@mid%
    \ifdim\lastskip<\insertskipamount\removelastskip\insertskip\fi
    \nointerlineskip\box\z@\penalty-200\insertskip
  \else%
    \SaveLastSkip
    \insert\topins{\penalty100 
    \splittopskip\z@skip
    \splitmaxdepth\maxdimen \floatingpenalty\z@
    \ifp@ge \dimen@\dp\z@
    \vbox to\vsize{\unvbox\z@\kern-\dimen@}
    \else \box\z@\nobreak\insertskip\fi}
    \RestoreLastSkip
   \fi\endgroup}
%
\catcode `\@=\newinsCatAt

\hoffset 0.5truein     
\voffset 1truein       
\def\ep{\varepsilon}
\def\inv{^{-1}}
\def\a{\longrightarrow}
\def\C{{\cal C}}
\def\pa{\partial}

\reflist

\refkey\ArMu {\bf Sh\^or\^o Araki}, {\bf Mitutaka Murayama}, 
{\it $G$-homotopy types of $G$-complexes and representations of
$G$-cohomology theories}, Publ. Res. Inst. Math. Sci. 14 (1978)
203--222

\refkey\Haus {\bf H.Hauschild}, {\it \"Aquvariante Konfigurationsr\"aume
und Abbildungsr\"aume}, Sprin\-ger Lecture Notes No. 788 (1979) 281--315

\refkey\KS  {\bf U.Koschorke, B.Sanderson}, {\it Self-intersections and
higher Hopf invariants}, Topology, 17 (1978) 283-290

\refkey\Dusa {\bf Dusa McDuff}, {\it Configuration spaces of
positive and negative particles}, Topology, 14 (1975) 91--107

\refkey\Ma 
{\bf J.P.May}, {\it The geometry of iterated loop spaces},
Springer Lectures Notes Series, {\bf 271} Springer--Verlag (1972)

\refkey\Mi 
{\bf R.J.Milgram}, {\it Iterated loop spaces}, Annals of Math., 
{\bf 84} (1966) 386--403

\refkey\RS {\bf Colin Rourke, Brian Sanderson}, {\it The
compression theorem},\nl {\tt arxiv:math.GT/9712235}

\refkey\SInv {\bf Graeme Segal}, {\it Configuration spaces and   
iterated loop spaces}, Invent. Math. 21 (1973) 213--221

\refkey\SPre {\bf Graeme Segal}, {\it Some results in equivariant
homotopy theory}, preprint (1979) available from:
{\tt http://www.maths.warwick.ac.uk/\char'176bjs/segal.html}

\refkey\Waner {\bf Stefan Waner}, {\it Equivariant homotopy theory 
and Milnor's theorem}, Trans. Amer. Math. Soc. 258 (1980) 351--368

\endreflist

\title{Equivariant configuration spaces}
\author{Colin Rourke\\Brian Sanderson}
\address{Mathematics Institute\\University of Warwick\\
Coventry\ \ CV4 7AL\ \ UK}
\email{cpr@maths.warwick.ac.uk\\ bjs@maths.warwick.ac.uk}

\abstract 
We use the compression theorem (cf [\RS; section 7]) to prove results
for equivariant configuration spaces analogous to the well-known
non-equiv\-ari\-ant results of May, Milgram and Segal [\Ma,\Mi,\SInv]. 
\endabstract

\primaryclass{55P91, 55P35, 55P40}

\secondaryclass{57R91, 55P45, 55P47}

\keywords{Equivariant homotopy theory, configuration spaces, 
multiple loop spaces, multiple suspension, group completion}  

\makeshorttitle

This paper is concerned with the homotopy type of equivariant
configuration spaces.  Some of the results proved here are known or
implicit in known results [\Haus,\Dusa,\SInv,\SPre].  However there
does not appear to be a systematic treatment in full generality in the
literature.  Moreover the treatment given here is new and several of
our results are deduced from the compression theorem [\RS] by
short, elementary arguments.  (The proof of the compression theorem in
[\RS] is also short and elementary.)

We are grateful to Peter May for pointing out this gap in the
literature and suggesting that the compression theorem might be used
to fill it and we are grateful to Graeme Segal for a copy of his
preprint [\SPre].  We would also like to thank the referee for many
helpful comments which have resulted in a clearer exposition of our
material.

Throughout $G$ is a finite group and $V$ is a finite-dimensional 
real representation of $G$.  

\section{Notation and results}

A {\sl based} $G$--space $X$ is a $G$--space with basepoint $*$ which is
fixed under the action.  We shall assume that all base points are 
non-degenerate.  For a $G$--space this means that the 
inclusion $\{*\}\subset X$ satisfies the HEP through equivariant maps 
(or equivalently that there is an equivariant retraction of $X\times I$
on $X\times \{0\}\cup \{*\}\times I$).  

Let $x\in X$, then $G_x$ denotes the {\sl isotropy group} of $x$ in
other words the subgroup $\{g\in G\mid g(x)=x\}$ of $G$.  Let $H$ be any
subgroup of $G$, then $X^H:=\{x\in X\mid G_x\supset H\}$ (the
fixed-point set of $H$) and $X_H:=\{x\in X\mid G_x=H\}$.  We say that
$X$ is $G$--{\sl connected} if $X^H$ is connected for each subgroup
$H$ of $G$.

We denote by $S^V$ the $G$--manifold which is the 1-point compactification
of $V$ (thought of as Euclidean space with $G$ action) with basepoint
at $\infty$.  If $X$ is a based $G$--space then $\Omega^V(X)$ is the
space of based $G$--equivariant maps $S^V\to X$ and $S^V(X)$ is the
$G$--space $S^V\wedge X$ with product action.

Let $M$ be a $G$--manifold without boundary.  Let $C_M(X)$ denote the
space of $G$--equivariant configurations in $M$ with labels in a based
$G$--space $X$.  Thus a point of $C_M(X)$ is a finite subset of $M$
which is closed under the action of $G$ and is equivariantly labelled
in $X$.  The basepoint of $C_M(X)$ is the empty configuration and the
topology is such that points labelled by $*$ can be deleted.  More
precisely, a labelled configuration of $q$ points is a point of the
quotient of the $q$--fold cartesian product
{\large$\times$}$\!_q(M\times X)$ by the symmetric group $\Sigma_q$,
acting by diagonal action, such that the $q$ points of $M$ are
distinct.  The configuration is equivariant if whenever $(m,x)$ is one
of the labelled points then so is $(g(m),g(x))$ for all $g\in G$.  We
identify configurations which agree apart from points labelled by $*$.
The topology on $C_M(X)$ is the quotient topology.  In other words
$C_M(X)$ is a quotient space of a subset of
$\coprod_q(${\large$\times$}$\!_q (M\times X))$.

There is a well-known map $j_V\co C_V(X)\to \Omega^V S^V(X)$ which can
be described by thinking of the configuration as a set of charged
particles and using the electric field they generate (see [\SInv; page
213]).  There is a more useful description of $j_V$ given by replacing
$C_M(X)$ by the space $C_M^\circ(X)$ of ``little discs'' in $M$
labelled in $X$. To be precise, suppose that $M$ is equipped with an
equivariant metric.  A {\sl little disc} at $x\in M$ is an $\ep$--disc
$D_\ep$ with centre $x$ where $\ep$ is small enough that the
exponential map is a diffeomorphism from the $\ep$--disc in the fibre
of $T(M)$ at $x$ to $D_\ep$. A point of the little disc space
$C_M^\circ(X)$ is a finite disjoint collection of little discs in $M$
which is closed under the action of $G$ and such that each little disc
is labelled by a point of $X$ and this labelling is equivariant.  The
topology is such that discs labelled by $*$ are ignored (as above).
The map $C_M^\circ(X)\to C_M(X)$ which takes each disc to its centre
is well-known to be a homotopy equivalence, see eg [\KS; Appendix on
page 289].

Replacing $C_M(X)$ by $C_M^\circ(X)$, the map $j_V$ can be described 
as follows.  A little
disc at $x\in V$ can be canonically identified with the unit disc $D$ in
$V$ (via translation and radial expansion).
Given a collection of little discs in $V$ define a map $S^V\to S^V$
by mapping each little disc to $S^V$ by using the identification
with $D$ and the standard identification $D/\pa D = S^V$.  Map the
rest of $V$ to $*$.

\proclaim{Theorem 1}Let $X$ be a based $G$--connected $G$--space.  
Then $j_V\co C_V(X)\to \Omega^V S^V(X)$ is a weak homotopy equivalence.\rm
\ppar

Notice that in theorem 1 neither $C_V(X)$ nor $\Omega^V S^V(X)$ is
necessarily a monoid.  However if we let $V+1$ denote the
representation $V\oplus\hbox{id}_\re$ then $C_{V+1}(X)$ is equivalent
to a monoid and we can form the classifying space $BC_{V+1}(X)$.

\proclaim{Theorem 2}Let $X$ be a based $G$--space. Then $j_{V+1}$ 
induces a weak homotopy equivalence $BC_{V+1}(X)\to \Omega^V 
S^{V+1}(X)$.\rm
\ppar

Note that in theorem 2 and the following corollary $X$ is not assumed
to be $G$--connected.  Recall that if $A$ is a topological monoid then
the natural map $A\to\Omega BA$ is called the {\sl group completion}
of $A$.  Group completion is a homotopy notion; thus we also refer to
a group completion composed with a homotopy equivalence as a group
completion.  A group completion composed with a {\it weak} homotopy
equivalence is called a {\sl weak group completion} (cf remark 1 at
the end of this section).

\proclaim{Corollary}Let $W$ be a representation with a trivial summand.
Then \break $j_W\co C_W(X)\to \Omega^W S^W(X)$ is a weak group completion.

\prf  Write $W=V+1$ and apply theorem 2. \qed

\ppar
Theorem 1 is new.  The special case when $V$ has a trivial
summand and the corollary to theorem 2 are proved in Hauschild [\Haus]
(Hauschild proves the result for $X=S^0$, however his methods extend 
with little alteration, see the remarks at the end of [\Haus]).
A generalisation of theorem 1 (where $V$ is replaced by
an arbitrary $G$--manifold) is given in section 4.
Theorem 2 can be deduced using Hauschild's methods [\Haus], but we prefer
instead to deduce it using Segal's methods [\SInv], which can be adapted to
prove the following stronger result:

\proclaim{Theorem 3}Let $M$ be a $G$--manifold.  Then
$BC_{M\times\re}(X)$ is homotopy equivalent to $C_M(SX)$.\rm
\ppar

In the case that $X$ is not necessarily $G$--connected and $V$ has no
trivial summand, then the above results give us no information about
$C_V(X)$. However in this case we can write $C_V(X)$ as $X^G\times A$
where $A$ is a certain topological monoid (for details here see
section 4).  We then have:

\proclaim{Theorem 4}There is a weak homotopy equivalence $t\co\Omega^VS^V(X)
\to X^G\times \Omega BA$ such that $t\circ j_V$ is homotopic to the
natural map $C_V(X)=X^G\times A\to X^G\times \Omega BA$.\rm \ppar 

This result improves on Segal's result [\SPre; theorem B].  Segal
proves a stable result for the special case  $X=S^0$.  Theorem 4 
is significantly sharper even in this special case.

\rk{Remarks\qua1}We have been careful to state exactly when we can 
prove homotopy equivalence and when we can only prove weak homotopy
equivalence.  It is worth noting that if $X$ has the $G$--homotopy
type of a $G$--CW complex then all the spaces defined above have the
homotopy type of CW complexes; in this case the word ``weak'' can be
deleted from all the results (cf Araki--Murayama [\ArMu] and Waner
[\Waner]).

{\bf2}\qua As remarked below theorem 1, equivariant configuration
spaces and loops--suspension spaces are not usually monoids.  This
contrasts with the situation in the classical loops--suspension
theorem where the configuration space approximates the monoid
structure closely (see eg May [\Ma]).  The one case where the
configuration space $C_V(X)$ is a monoid (up to homotopy) is when the
representation $V$ has a trivial summand.  In this case
$\Omega^VS^V(X)=\Omega(\Omega^{V'}S^V(X))$ is a loop space (where
$V=V'+1$).  The equivalent monoid structure on $C_V(X)$ is described
in detail in section 3 (see the start of the description of $q$) and
this description makes it clear that $j_V$ commutes (up to homotopy)
with the monoid structures.  This fact is implicit in theorem 2.

\section{Proof of theorem 1}

We follow the proof given in [\RS] of the non-equivariant case.  Here
is a quick sketch of the argument.  A representative of $\pi_kC_V(X)$
is a parallel framed manifold in $\re^k\times V$ equivariantly labelled
in $X$, whilst a representative of $\pi_k\Omega^V S^V(X)$ is a
semi-parallel framed manifold.  The compression theorem applied
locally implies that they are equivalent.  The description of
$\pi_k\Omega^V S^V(X)$ is justified by transversality to $0\in V$.
The concepts ``parallel'' and ``semi-parallel'' are defined in the
detailed version of the argument, which now follows.  We start with
the details of the transversality result that we need.

Consider a $G$--map $f\co V\times Q\to V\times X$, where $Q$ is a
smooth manifold with trivial $G$--action, such that the composition
with projection on $V$ is smooth.  Let $x\in V\times Q$ and let
$H=G_x$ (the isotropy group of $x$) and let $W=V^H$ (the subspace of
$V$ fixed by $H$).  So $x\in W\times Q$.  Let $U$ be the orthogonal
complement of $W$ in $V$ (we assume without loss that we have an
equivariant Riemannian metric on all $G$--manifolds) then $U$ can be
identified with each fibre of the normal bundle of $W\times Q$ in
$V\times Q$.  Denote the fibre at $x$ by $U_x$.

Now $f(x)\in W\times X$ and we
can again identify $U$ with each fibre of the normal bundle
of $W\times X$ in $V\times X$.  Denote the fibre at $f(x)$ by
$U_{fx}$.  We say that $f$ is {\sl semi-parallel} at $x$ if there
is a neighbourhood $N$ of $0$ in $U$ such that $f\vert N_x$ is
the identity $N_x\to N_{fx}$, where $N_x,N_{fx}$ are the copies of
$N$ in $U_x,U_{fx}$ respectively.  There is a similar definition of
semi-parallel near a subset of $V_H\times Q$ which is closed in
$W\times Q$.  (Recall that $V_H\subset W=V^H$ is the subset of points
with isotropy group $H$.)

We say that $f$ is {\sl transverse} to $0\times X\subset V\times X$ if
for each subgroup $H$ of $G$ we have that $f\vert V_H\times Q$ is
transverse to $0\times X\subset V\times X$ and further $f\inv(0\times X)
\cap V_H\times Q$ is a closed subset $C$ of $V^H\times Q$ (in fact a
submanifold by transversality) and $f$ is semi-parallel near $C$.
(Notice that the definition implies that $f\inv(0\times X)$ is a 
disjoint collection of manifolds each of which has the property
that all points have the same isotropy group.)

\proclaim{Transversality lemma}The map $f$ is $\ep$--equivariantly
homotopic to a transverse map.  Further if $f$ is already transverse
near a closed subset $C$ (invariant under $G$) then we can assume the 
homotopy fixed on $C$.

\prf  Consider the conjugacy classes $\{G\}=\C_1,\C_2,\ldots,
\C_q=\{\{e\}\}$ of subgroups of $G$ ordered such that $H\in \C_i$,
$K\in \C_j$, $H\subset K$ imply $j<i$.
Suppose that $f$ is transverse near $f\inv(0\times X)\cap V^{H}$ 
for each $H\in \C_i$, $1\le i<j$ and consider $H\in \C=\C_j$. Let $W=V^H$ 
and $U$ be as above and consider the open subset $V_H\subset W$.  
Denote by $V_\C$ the union $\cup_{H\in\C}V_H$ and by $W_\C$ the 
union $\cup_{H\in\C}V^H$. 
By supposition $f$ is already transverse near $W-V_H$ for each
$H\in\C$ and this implies that $f\inv(0\times X)
\cap V_\C$ is a closed subset $C$ of $W_\C$.  

We start by making $f$ semi-parallel near $C$.  For any point 
$x\in V^H$, $H\in\C$ we can deform $f$ to be semi-parallel by choosing 
a small $\ep$--neighbourhood
$N$ of $0\in U$, squeezing $f\vert N_x$ to $f(x)$ and growing
$f\vert N$ radially outwards in $U_{fx}$ to be the identity 
$N_x\to N_{fx}$.  By doing the same homotopy for each point
in the orbit of $x$ (which is a discrete subset of $V_\C$ in
bijection with the set of cosets $G/H$) we have an equivariant
homotopy which makes $f$ semi-parallel at $x$.  A similar
process will make $f$ semi-parallel near $x$ (work in a 
neighbourhood which is small enough not to meet any of its
translates under $G/H$) and we can leave $f$ fixed near a 
closed subset where it is already semi-parallel.  A standard
patch-by-patch argument will now make $f$ semi-parallel near $C$.

We now apply ordinary transversality to $f\co V_\C\times Q\to W_\C\times X$
to make it transverse to $0\times X$ by a homotopy which leaves the
$X$--coordinate fixed.  Again we work patch by
patch, copying the (small) homotopy around orbits, and preserve
the semi-parallel condition by extending the movement by crossing
with the identity on $N$ for a suitably small neighbourhood of
$0$ in $U$.  The method of proof makes the relative statement
clear. \qed\ppar

Continuing with the proof of theorem 1, we observe that a 
representative of $\pi_kC_V(X)$ is a collection $T$ of parallel 
framed manifolds in $S^k\times W^H$ equivariantly labelled in $X$
for varying $H$, where {\sl parallel} means having normal bundle
whose fibres can be identified with the appropriate translate of
$V$.  Further, we may assume that these are smooth manifolds with
boundary labelled by $*$.  This uses non-degeneracy of the basepoint,
which implies that $X$ is homotopy equivalent to $X$ with a whisker 
$\{*\}\times I$ attached at the basepoint.  Deform the labelling
map $l$ to be transverse to $*\times 1/2$ and then stretch $[0,1/2]$
to $[0,1]$; the proof of transversality needed here is a simplified 
version of the proof given above.  Make $l\vert T^{H_i}$ transverse
to $*\times 1/2$ by induction on $i$;  each move is extended to
a neighbourhood by moving normal fibres like their centres and the
process then continues.

Now by the main transversality lemma, a representative
of $\pi_k\Omega^V S^V(X)$ is a collection $U$ of semi-parallel framed 
manifolds.  Notice that the argument just given implies that we
may ignore the labelling map near $*$ and hence that we are
mapping into $V\times X$ rather than $S^V(X)$.  Moreover again by
the same argument we can assume that these manifolds have boundaries 
labelled by $*\times1/2$, which is then replaced by $*$.  But 
$G$--connectivity of $X$ implies that all components of these 
manifolds can be assumed to have boundary.  Choose a point $x\in U_H$
say and choose a small disc centre $x$, with collar both invariant
under $H$ where small means that the disc does not meet any of 
its translates under representatives of $G/H$.  Use connectivity of
$X^H$ to deform the label on the disc to $*$ and extend using the
collar and a path to the basepoint.  Do the same for the translates.
The interiors of the discs are now labelled by $*$ and can be
deleted.

Now by the compression theorem such a framed manifold (with boundary) 
can be made parallel; the version of the compression theorem needed here is
the codimension 0 (with boundary) version of the multi-compression
theorem (see the note at the end of section 4 of [\RS]).  The 
isotopy is not small in this version but of the form,
shrink to the neighbourhood of a codimension 1 spine, then perform a
small isotopy.  This can be made equivariant by choosing an
appropriate spine in the quotient manifold under the group action,
using an equivariant shrink and then the above patch by patch
argument.  \qed

\section{Proofs of theorems 2 and 3}

We shall describe a map $q\co BC_{M\times\re}(X)\to C_M(SX)$ which we
shall prove to be a homotopy equivalence; this proves theorem 3.
The description of $q$ will make it clear that the following diagram
commutes up to homotopy:
$$
\matrix{
C_{M\times\re}(X)&\buildrel j\over\a&\Omega C_M(SX)\cr
\cr           
\downarrow&&\vert\vert\cr
\cr
\Omega BC_{M\times\re}(X)&\buildrel\Omega q\over\a&\Omega C_M(SX)\cr
}
$$ 
Here $j$ is defined in a similar way to $j_V$; we replace
$C_{M\times\re}(X)$ by the homotopy equivalent space
$C_{M\times\re}^-(X)$ of disjoint finite sets of ``little 1-disks''
parallel to $\re$ (little bi-collars) and use a similar formula to
that for $j_V$ (above the statement of theorem 1).  Below we give a
precise definition of $BC_{M\times\re}$ in which $C_{M\times\re}$ is
replaced by a strictly associative monoid $C$.  The vertical
unlabelled map in the diagram is then equivalent to the natural map
$C\to \Omega BC$.  Thus the homotopy equivalence $q$ is induced by $j$
and theorem 2 can be deduced by first applying theorem 3 with $M=V$
and then applying theorem 1 with $X$ replaced by $SX$.

\rk{Description of $q$}

In order to define the classifying space $BC_{M\times\re}(X)$ we shall
replace $C_{M\times\re}(X)$ by a homotopy equivalent strictly
associative mon\-oid and we do this in a standard way, cf [\SInv;
bottom of page 213].  We shall also use ``little 1-disks'' as above.
To be precise, consider the associative monoid $C$ of which a point is
a collection (configuration) of little 1-discs in $M\times (0,t)$ for
some $t\ge0$.  The 1-discs are required to be parallel to $\re$ and
the collection is closed under action of $G$ and equivariantly
labelled in $X$.  Addition is defined by juxtaposing the intervals
$(0,t), (0,t')$ to form the interval $(0,t+t')$.

The classifying space $B=BC$ is the realisation of the simplicial
space of which an $n$--simplex is an $n$--tuple of such
configurations, with faces defined by omitting the first or last
configuration or by adding an adjacent pair of configurations.
We shall define $q\co B\to C_M(SX)$.  

First we require some technicalities on slices of the standard simplex.\par
Let $\Delta_n = \{\{ x_1, \dots, x_n\} \in \re^n\mid  0\le x_n\le\cdots
\le x_1\le 1\}$ denote the standard $n$-simplex.  
The {\sl $k^{th}$ slice of $\Delta_n$}, denoted $L_k$ is the $(n-1)$--cell
$\{x\in\Delta_n \mid x_k={1\over 2}\}$.  The
vertices of $\Delta_n$ are $v_0=(0,\ldots,0),v_1=(1,0,\ldots,0),
\cdots v_n=(1,\ldots,1)$.  Let $\sigma_k$ be the face of $\Delta_n$
spanned by $v_0,\cdots,v_{k-1}$ and let $\tau_k$ be the face spanned by
$v_{k},\cdots,v_n$.  Then $\Delta_n$ is the join of $\sigma_k$ with
$\tau_k$ and the product of $\sigma_k$ with $\tau_k$ can be identified
with $L_k$ by $(u,v)\mapsto {1\over 2}(u+v)$.  Given parameters
$s_k\in (0,\infty)$ for $1\leq k\leq n$ we will show how to thicken
and pull apart the slices $L_k$ so that interiors of thickened slices
are pairwise disjoint, see the figure.
The {\sl thickened $k^{th}$ slice}
$TL_k$ will be defined as the image of an embedding $S_k\co
(0,s_k)\times\sigma_k\times \tau_k\to \Delta_n$.

\fig{: disjoint thickened slices of the standard 3--simplex}
\relabelbox\small
\epsfysize 2.2truein
\epsfbox{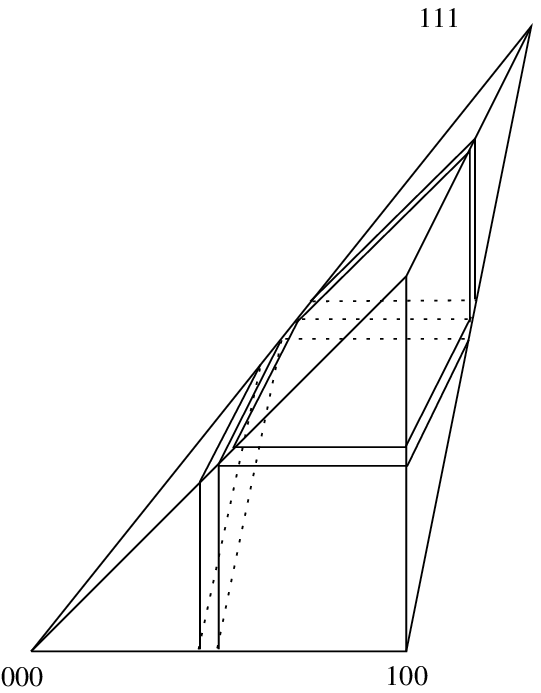}
\adjustrelabel <-7pt, 0pt> {000}{$\scriptstyle(0,0,0)$} 
\adjustrelabel <-7pt, 0pt> {100}{$\scriptstyle(1,0,0)$} 
\adjustrelabel <12pt, 3pt> {111}{$\scriptstyle(1,1,1)$} 
\extralabel <-28pt, 1.19in> 
{$\lower 1pt\hbox to 20pt{\leftarrowfill}\scriptstyle(1,1,0)$} 
\endrelabelbox
\endfig

Let $\theta\co [0,\infty) \to [{1\over 2},1)$ be the homeomorphism
 given by $\theta(t)={2\over\pi}\arctan(1+t)$.  Let $v_i \in \sigma_k,
 v_j \in \tau_k$, so $i < k\le j$.  Define $S_k(t,v_i ,v_j)=u_{i,k}v_i
 +(1-u_{i ,k})v_j$, where $u_{i,k}=\theta(s_{i +1}+\cdots
 +s_{k-1}+t)$, and in general define
$$S_k(t,u,v)=
\sum_{\scriptstyle0\leq i< k\atop\scriptstyle k\le j\le n}\lambda_i \mu_j S_k(t,v_i ,v_j)$$
where $u=\sum_{0\le i < k}\lambda_i v_i$ and $v=\sum_{k\le j\le n} \mu_jv_j$.

We are ready to define $q$.  Consider configurations $c_k$ in $M\times
(0,s_k)$, for $1\le k\le n$, and labelling functions $x_k\co c_k\to
X$. Each little 1-disc $d\in c_k$ is determined by its centre point
$c(d)$ and the radius $\ep(d)$ of the 1-disc.  A typical
point in an $n$-simplex of $B$ can be written as
$b=[t_1,\cdots,t_n,(c_1,x_1),\cdots,(c_n,x_n)].$ We will define
$q(b)$. If $(t_1,\cdots,t_n)$ is not in the image of any $S_k$ then
define $q(b)=*$.  Suppose then that $(t_1,\cdots,t_n)=S_k(t,u,v)$.
Let $d\in c_k$, then $c(d)=(c',c'')\in M\times (0,s_k)$ and we define
$c'$ to be in the configuration of $q(b)$ if $|t-c''|\le \ep(d)$. The
label of $c'$ in $SX$ has $X$-coordinate $x_k(d)$ and suspension
coordinate ${1\over 2\ep(d)} (t+\ep(d) -c'')$.

\rk{Proof that $q$ is a homotopy equivalence}

We adapt the proof given in Segal [\SInv] of the case when $M=\re^n$
and the action of $G$ is trivial (see [\SInv; proposition 2.1]).  We need to
consider two other spaces.  $B'$ is the classifying space of the
partial monoid of configurations in $M$ labelled in $X$, where the
operation is union of disjoint configurations (see [\SInv; page 215])
and $B''$ is the classifying space of the sub-partial monoid of
$C\simeq C_{M\times\re}(X)$ (defined above) comprising configurations
which project homeomorphically to $M$ with the operation being defined
only for configurations which project disjointly to $M$.  Projection
$p\co B''\to B'$ is a homotopy equivalence with inverse induced by
inclusion of $M$ as $M\times {1\over2}\subset M\times(0,1)$.  Further
$q'\co B'\to C_M(SX)$ can be defined by using the halfway slices
$L_k$ and ``full-size'' bicollars
and then $q'\co B'\to C_M(SX)$ is a homeomorphism by the same proof as
[\SInv; proposition 2.3].  Thus it suffices to prove that the
inclusion $i\co B'' \to B$ is a homotopy equivalence, and that the
following diagram commutes up to homotopy:
$$\matrix{
B''&\buildrel p\over\a& B'\cr
\cr           
\downarrow i&&\downarrow q'\cr
\cr
B&\buildrel q\over\a& C_M(SX)\cr
}$$
The following considerations explain why it is to be expected that
$i\co B'' \to B$ is a homotopy equivalence.

A thin $\re$--slice of a configuration in $M\times\re$ projects
homeomorphically to $M$.  Further the addition in $C$
shows that any configuration is the sum of such configurations.
This implies that any 1--simplex in $B$ can be deformed across
a sequence of 2--simplexes until it is replaced by a chain
of 1--simplexes in $B''$.  In a similar way, higher-dimensional
simplexes in $B$ can be deformed into $B''$.  

This argument can readily be sharpened to prove that $B''\to B$ 
is a weak homotopy equivalence.  The proof given in [\SInv] can be
regarded as sharpening this argument to prove a full homotopy
equivalence.  In [\SInv] the process of splitting configurations into
$\re$--slices is made precise as ``disintegration'', see [\SInv;  
pages 217--218].
On a more formal level, the proof in [\SInv] carries over to the
equivariant case with obvious changes.  There is one comment which
needs to be made.  We need an equivariant version of the application
of non-degeneracy of the basepoint made just below [\SInv, proposition 2.7].
This is done by using a simple transversality argument as in the
proof of theorem 1 above. 

Finally to see that $q'p\simeq qi$ observe that if $b\in B''$ is
determined as above by an $n$-tuple of configurations of little
1-discs then $q'p(b)$ is defined using the halfway slices and
full-size bicollars whilst $qi(b)=q(b)$, defined above, uses slices
near but not exactly halfway and with small bicollars.  Thus the
required homotopy is given by sliding the slices to halfway (use
the formula for the slices above and a homotopy of $\theta$ to the
constant map to $1\over2$) and expanding the collars.  \qed

\section{Proof of theorem 4}

We shall need an extension of theorem 1 which applies to configurations 
in a general $G$--manifold.

Let $M$ be a $G$--manifold of dimension $n$ without boundary and with
tangent bundle $T_M$.  Let $E_M$ be the sphere bundle in $T(M)\oplus
1$; this bundle can also be thought of as the disc bundle in $T(M)$
with each fibre of the sphere bundle squeezed to a point (which is
taken as the basepoint in the fibre).  Let $E_M(X)$ denote the bundle
with fibre $S^n(X)$ obtained from $E_M$ by taking the smash product
with $X$ in each fibre.  Let $\Gamma_M(X)$ denote the space of
equivariant sections of $E_M(X)$.  $\Gamma_M(X)$ is based by the base
section (given by the basepoint in each fibre).

There is a map $j\co C_M(X)\to\Gamma_M(X)$ given by replacing $C_M(X)$
by the equivalent little disc space $C_M^\circ(X)$ and then using
the description given in 
McDuff [\Dusa; page 95].  (McDuff calls this map $\phi_\ep$; the
formula is in the middle of the page.)  Theorem 1 has the following 
generalisation:

\proclaim{Theorem 5}Let $X$ be a based $G$--connected space.  Then $j\co
C_M(X)\to\Gamma_M(X)$ is a weak homotopy equivalence.\rm

\prf
The proof is a straightforward generalisation of the proof of theorem
1.  We first discuss the case when $G$ is trivial.  An element of
$\pi_kC_M(X)$ can be regarded as a compact submanifold of $\re^k\times
M$ of dimension $k$, such that the projection on $\re^k$ is an
immersion, and which is labelled in $X$ so that the boundary is
labelled by $*$ (see the whisker argument used in the proof of theorem
1 here).

Now an element of $\pi_k\Gamma_M(X)$ can be interpreted by transversality
to the zero section as a compact submanifold $W$ say of $\re^k\times
M$ which is labelled in $X$ so that the boundary is labelled by $*$.
Further the normal bundle pulls back from $T(M)$ and hence each fibre
has an explicit identification with the tangent fibre to $M$ at
the same point.  Locally, using a local trivialisation of $T(M)$, this 
means that the normal bundle is framed with framing not necessarily
the same as this local trivialisation.  But the compression theorem
can now be applied locally to make these framings agree and thus $W$
can locally be moved to a parallel submanifold.  Using the relative
version of the compression theorem and working patch by patch we can move 
$W$ globally to a parallel submanifold and this shows that $j_*\co
\pi_kC_M(X)\to \pi_k\Gamma_M(X)$ is surjective.  A similar relative
argument shows that $j_*$ is injective.

This proof is extended to the general case (when $G$ is non-trivial)
in exactly the same way that the proof of the non-equivariant
loops--suspension theorem given in [\RS; section 7] was extended to
the equivariant case in the proof of theorem 1:  we replace
transversality by equivariant transversality and we replace the local
moves given by the compression theorem by equivariant moves, by copying
the local move around an orbit. \qed

\ppar
We can now prove theorem 4.  If $V$ has no non-trivial summand,
we can write  $C_V(X)$ as the product $X^G\times C_{Q\times\re}(X)$
where $Q$ is the unit sphere in $V$.  This is because 
$0\in V$ must
be labelled by a fixed point of $G$ in $X$ and the configuration is
described by giving the label at 0 (which might be $*$) and an 
equivariant configuration in $V-\{0\}\cong Q\times\re$.
Write $A$ for $C_{Q\times\re}(X)$

Now by theorem 3, there is a group completion $A\to
\Omega C_Q(SX)$.  Further by theorem 5 $C_Q(SX)$ has the weak homotopy
type of $\Gamma_Q(SX)$.  Thus there is a weak group completion 
$A\to\Omega\Gamma_Q(SX)$.

Next consider $\Gamma_Q(SX)$.   This is the space of cross-sections of
$E_Q(SX)$, where $E_Q$ is the disk bundle of $T_Q$, with
each boundary sphere squeezed to a point, 
and $E_Q(SX)$ is obtained by smashing each fibre with 
$SX$.  Ignoring the $G$--action, this is a bundle with fibre 
$S^{n-1}\wedge S^1\wedge X$ and, ignoring $X$ as well for the moment,
we can identify each fibre $S^{n-1}\wedge S^1$ with $S^V$:  think
of the $S^{n-1}$ in the fibre over $q\in Q$ as the subspace of $V$
parallel to the tangent plane to $Q$ at $q$, compactified at $\infty$;
think of $S^1$ as the similarly
compactified perpendicular line.  It follows that we can identify
each fibre of $E_Q(SX)$ with $S^V(X)$ and taking account of
the $G$--action we see that 
$\Gamma_Q(SX)$ is the space of equivariant maps $Q\to S^V(X)$.

We can now interpret $\Omega\Gamma_Q(SX)$ as the subspace $Z$ of 
$\Omega^VS^V(X)$ comprising maps $S^V\to S^V(X)$ which near zero 
and infinity map to $*$.
We shall construct a homotopy equivalence $u\co \Omega^VS^V(X)\to 
X^G\times Z$.  The required weak homotopy equivalence $t$ (of the
statement of theorem 4) is the composition of $u$ with the identification 
of $Z$ (up to weak homotopy type) with the group completion of 
$A$ which is described above.  It is straightforward to check 
that $t\circ j_V$ is homotopic to the natural map 
$X^G\times A\to X^G\times Z$ and theorem 4 follows.

It remains to construct $u$.  Consider an equivariant map $f\co S^V\to
S^V(X)$.  Now $f$ must map $0\in V$ to either the basepoint of
$S^V(X)$ or to $0\times p$ in $S^V\wedge X$ where $p\in X^G$ and $p\ne
*$.  Thus in either case $f$ determines a point of $X^G$.  Further, by
squeezing a neighbourhood of $0$ to $0$ and then radially expanding a
small sphere near zero to infinity, we can homotope $f$ so that there
is a small disc $D$ centred at $0$ such that $f$ identifies $D/\pa D$
with $S^V\times p$.
Further we can squeeze a neighbourhood of $\infty$ to $\infty$.  Once
this has been done, we can see that $f$ on $V-D$ determines a map in
$Z$.  These constructions define $u\co \Omega^VS^V(X)\to X^G\times Z$.
There is a map $v\co X^G\times Z\to\Omega^VS^V(X)$ given by sending
$(p,g)$ to the map which identifies $D/\pa D$ with $S^V\times p$
and maps $V-D$ by $g$.  It is straightforward to check that $u$ and
$v$ are inverse homotopy equivalences. \qed\np

\references
\bye